\documentclass[10pt,a4paper,reqno]{amsart}
\usepackage{amssymb}
\usepackage{amsmath}
\usepackage{enumerate}
\usepackage{graphicx}
\usepackage{multirow}
\newfont{\jine}{wncyb10}
\newcommand{\re}{\operatorname{reg}}
\newcommand{\au}{\operatorname{Aut}}

\begin{document}

\title[On some connected groups of automorphisms of Weil algebras]{On some connected groups of automorphisms of Weil algebras}
\author{Miroslav Kure\v{s} and Jan \v{S}\' utora}
\address{
 Brno University of Technology \newline
\indent Department of Mathematics \newline
\indent Technick\'a 2, 61669, Brno, Czechia}
\email{kures@fme.vutbr.cz}
\address{
 Brno University of Technology \newline
\indent Department of Mathematics \newline
\indent Technick\'a 2, 61669, Brno, Czechia}
\email{sutora@post.cz}

\setcounter{page}{0}


\begin{abstract}
The method of direct calculation of the group of $\mathbb R$-algebra automorphisms of a Weil algebra is presented in detail. 
The paper is focused on the case of a one-componental group and presents two cases of values of the determinant of its linear part.
\end{abstract}

\maketitle

\section*{Introduction}
The paper is devoted to groups of $\mathbb R$-algebra automorphisms  of Weil algebras, which are important examples of algebraic groups.
We consider the usual Euclidean topology.
The essential fact (\cite{ASA}) is that every algebraic group has as irreducible components the cosets modulo a closed connected normal subgroup of
finite index, called the {\it identity component}. (It contains the neutral element which is the identity automorphism.)
Those irreducible components
turn out to be also the connected components. 
We denote the group of $\mathbb R$-algebra automorphisms of a Weil algebra $A$ by $\operatorname{Aut} A$ and its
identity component by $G_A$.

In differential geometry, Weil algebras play a fundamental role in the theory of Weil bundles $T^AM$ over a smooth manifold $M$
which generalize well-known higher order velocities bundles $T^r_nM$.
Contact elements are defined as orbits of the group action. 
In the classical case, they form the fiber bundle denoted by $K^1_nM=\re T^1_n M/G^1_n$
where $G^1_n$ denotes the first order jet group. 
It is the bundle of contact elements; moreover the bundle of oriented contact elements is obtained
if we replace the whole group $G^1_n$ by its identity component which is the subgroup of matrices with positive determinant, see \cite{GUG}.
We refer to paper 
\cite{KUR}
for the case of higher order jet groups $G^r_n$, bundles $K^r_nM=\re T^r_n M/G^r_n$ of higher order contact elements and the generalizations:
groups of Weil algebra $\mathbb R$-algebra automorphisms $\au A$ and bundles $K^AM=\re T^A M/\au A$ of Weil contact elements.

We investigate the case $\au A=G_A$, it means Weil algebras possessing only one connected component. They represent an important case
where orientation reversing reparametrizations
for corresponding velocities are impossible.

\section
{Basic concepts, denotations and the method}

This section will, among other things, show how to effectively calculate a group of $\mathbb R$-algebra automorphisms of a Weil algebra.

The Weil algebra $\mathbb D^2_2/\langle X^2,Y^2\rangle$ represents the Weil algebra of the iterated tangent functor. 
Let the residue classes of $X$ and $Y$ will be denoted by the same symbol $X$ and $Y$, respectively.
Then we choose the basis $\left\{ 1,X,Y,XY \right \}$ and write elements $a+bX+cY+dXY$ of this Weil algebra as vectors $(a,b,c,d)$.
The endomorphisms have a form
$$
(a,b,c,d)
\begin{pmatrix} 
1 & 0 & 0 & 0 \\
0 & A & B & C \\
0 & D & E & F \\
0 & G & H & I  \end{pmatrix}
=
(a,Ab+Dc+Gd,Bb+Ec+Hd,Cb+Fc+Id)
$$
which follows from the fact $1\mapsto 1$. It is therefore sufficient to investigate the third order submatrix 
$
\mathcal M=
\begin{pmatrix} 
A & B & C \\
D & E & F \\
G & H & I  \end{pmatrix}
$. 
In fact, we restrict to subspace formed by nilpotent elements (which is an ideal in the algebra and which is usually denoted by $\mathfrak n$).
For automorphisms, the matrix $\mathcal M$ must be non-singular. However, we have to describe {\sl algebra} endomorphisms and automorphisms. 
As
$$
X\mapsto AX+BY+CXY \;\;\text{ and }\;\;
Y\mapsto DX+EY+FXY,
$$
we compute
$$
XY\mapsto ADX^2+(AE+BD)XY+BEY^2
$$
and it follows $AD=0$ and $BE=0$. So, we have four variants
\begin{enumerate}
\item[(i)]$B=0$, $D=0$
\item[(ii)]$A=0$, $E=0$
\item[(iii)]$D=0$, $E=0$
\item[(iv)]$A=0$, $B=0$
\end{enumerate}
and, as $G=AD=0$, $H=BE=0$, $I=AE+BD$, they correspond with matrices
$$
\begin{pmatrix} 
A & 0 & C \\
0 & E & F \\
0 & 0 & AE  \end{pmatrix},
\begin{pmatrix} 
0 & B & C \\
D & 0 & F \\
0 & 0 & BD  \end{pmatrix},
\begin{pmatrix} 
A & B & C \\
0 & 0 & F \\
0 & 0 & 0  \end{pmatrix},
\begin{pmatrix} 
0 & 0 & C \\
D & E & F \\
0 & 0 & 0  \end{pmatrix}.
$$
It is clear that the third and fourth matrices will be killed in the case of automorphisms (they have zero row).
Thus, 
$\mathbb R$-algebra automorphisms are represented by matrices
$$
\begin{pmatrix} 
A & 0 & C \\
0 & E & F \\
0 & 0 & AE  \end{pmatrix},
\begin{pmatrix} 
0 & B & C \\
D & 0 & F \\
0 & 0 & BD  \end{pmatrix},
$$
in which $A\ne 0$, $E\ne 0$ and $B\ne 0$, $D\ne 0$, respectively. 
As non-zero coefficients can be either positive or negative, we find that the group
of $\mathbb R$-algebra automorphisms has eight connected components.
For the determinant $\mathfrak D=\det\mathcal M$,
we have
$$
\mathfrak D=A^2E^2 \quad\text{or}\quad \mathfrak D=-B^2D^2.
$$
If we make a deeper restriction, only for the supspace formed by elements of $\mathfrak n / \mathfrak n^2$,
we consider the matrix $\mathcal M_1$
having two possible forms 
$$
\begin{pmatrix} 
A & 0  \\
0 & E  \end{pmatrix},
\begin{pmatrix} 
0 & B  \\
D & 0   \end{pmatrix},
$$
($A\ne 0$, $E\ne 0$, $B\ne 0$, $D\ne 0$), and,
for the determinant $\mathfrak D_1=\det\mathcal M_1$,
we have
$$
\mathfrak D_1=AE \quad\text{or}\quad \mathfrak D_1=-BD.
$$

\section{Theorem on a connected group of $\mathbb R$-algebra automorphisms of a Weil algebra}

In the introductory example, we saw the group of automorphisms non-connected, $G_A$ was a proper
subgroup of $\au A$. We now want to focus only on algebras where the group is connected ($G_A=\au A$). 
We recall that in \cite{IVK} is proved that the Weil algebra
$\mathbb D^6_2/\langle X^3+Y^4,X^4+Y^5\rangle$
possesses
a connected (one-componental) group
of $\mathbb R$-algebra automorphisms
and
$$
\mathcal M_1=
\begin{pmatrix} 
1 & 0  \\
0 & 1  \end{pmatrix}
$$
(see Theorem 1 and its proof in \cite{IVK} for more details).
and
$$
\mathfrak D_1=1.
$$

The following assertion holds.

{{\sc Theorem.\,}\it
Let $A$ is a Weil algebra with connected group of $\mathbb R$-algebra automorphisms.
Then there may occur one of the following two cases for the determinant $\mathfrak D_1$:
\begin{enumerate}
\item[(i)] $\mathfrak D_1$ is constant with the value 1
\item[(ii)] $\mathfrak D_1$ takes all real values from the interval $(0,\infty)$.
\end{enumerate}}
\begin{proof}
It is obvious that the map assigning the determinant $\mathfrak D_1$
to an $\mathbb R$-algebra automorphism
is a continuous function and, moreover, a group homomorphism.
It is also clear that its image is a subset of the interval $(0,\infty)$.
So it is clear that only possibilities
(i) and (ii) can occur because we have no other connected subgroup
of the group $(0,\infty)$ (with the usual multiplication of reals) than these two.
The case (i) occurs for the Weil algebra
$$
\mathbb D^6_2/\langle X^3+Y^4,X^4+Y^5\rangle.
$$ 
It remains to show
that the case (ii) may occur, too.

Let us consider $A=\mathbb D^4_2 / \langle X^3Y, X^2Y^2, Y^4, X^3-Y^3 \rangle$. 
We choose the basis $$\left\{ 1,X,Y,X^2,XY,Y^2,X^3,X^2Y,XY^2,X^4 \right \}.$$
Then the endomorphism $\varphi\colon A\to A$
maps
$$
\phi(X)=AX+BY+CX^2+DXY+EY^2+FX^3+GX^2Y+HXY^2+IX^4
$$
and
$$
\phi(Y)=JX+KY+LX^2+MXY+NY^2+PX^3+QX^2Y+RXY^2+SX^4.
$$
For $\varphi$ to be an $\mathbb R$-algebra automorphism, it is necessary
\begin{enumerate}
\item[(i)]
$AK-BJ \ne 0$
\item[(ii)]
$\varphi(X^3Y)=\varphi(0)=0$, but it means that in the expression of $\varphi(X^3Y)$ the coefficient standing at $X^4-XY^3$ is zero, which gives the equation
\begin{equation}
A^3J-B^3J-3AB^2K=0
\label{e1}
\end{equation}
\item[(iii)]
$\varphi(X^2Y^2)=\varphi(0)=0$, but it means that in the expression of $\varphi(X^2Y^2)$ the coefficient standing at $X^4-XY^3$ is zero, which gives the equation
\begin{equation}
A^2J^2-2B^2JK-2ABK^2=0
\label{e2}
\end{equation}
\item[(iv)]
$\varphi(Y^4)=\varphi(0)=0$, but it means that in the expression of $\varphi(Y^4)$ the coefficient standing at $X^4-XY^3$ is zero, which gives the equation
\begin{equation}
J^4-4JK^3=0
\label{e3}
\end{equation}
\end{enumerate}
We start with the equation (\ref{e3}).
\begin{enumerate}
\item[($\alpha$)]
If $J=0$, then $A \ne 0$, $K \ne 0$ and then $B=0$ due to (\ref{e1}) and (\ref{e2}).
\item[($\beta$)]
If $J\ne 0$, then $J=\sqrt[3]{4} K$, (\ref{e1}) transforms to
\begin{equation}
\left(\sqrt[3]{4} A^3 - \sqrt[3]{4} B^3 - 3AB^2  \right)K=0
\label{e4}
\end{equation}
and 
(\ref{e2}) transforms to
\begin{equation}
\left(\sqrt[3]{8} A^2 - 2\sqrt[3]{4} B^2 - 2AB  \right)K^2=0,
\label{e5}
\end{equation}
but $K\ne 0$ as $J\ne 0$. So the expressions in brackets have to be zero
but the obtained system in unknowns $A$ and $B$ has only the solution $A=B=0$ which is impossible.
\end{enumerate}
So we continue with $B=0$, $J=0$, $A \ne 0$, $K \ne 0$.
Moreover, it is necessary
\begin{enumerate}
\item[(v)]
$\varphi(X^3-Y^3)=\varphi(0)=0$, but it means that in the expression of $\varphi(X^3-Y^3)$ the coefficient standing at $X^3-Y^3$ is zero, which gives the equation
\begin{equation}
A^3-K^3=0
\label{e6}
\end{equation}
with the solution $K=A$ and also in the expression of $\varphi(X^3-Y^3)$ the coefficient standing at $X^4-XY^3$ is zero, which gives the equation
\begin{equation}
3A^2CD-3A^2M=0
\label{e7}
\end{equation}
with the solution $M=C$.
\end{enumerate}
So, we have
\begin{eqnarray*}
\phi(X)&=&AX+CX^2+DXY+EY^2+FX^3+GX^2Y+HXY^2+IX^4 \\
\phi(Y)&=&AY+LX^2+CXY+NY^2+PX^3+QX^2Y+RXY^2+SX^4 \quad \text{where $A\ne 0$}.
\end{eqnarray*}
That is why we obtain
$$
\mathcal M=
\begin{pmatrix} 
A      & 0      & C        & D        & E         & F       & G      & H     & I   \\
0      & A      & L        & C        & N         & P       & Q      & R     & S   \\
0      & 0      & A^2      & 0        & 0         & 2AC     & 2AD    & 2AE   & 2AF+C^2-2DE \\
0      & 0      & 0        & A^2      & 0         & AL-AE   & 2AC    & AD+AN & AP+CL-AH-CE-DN \\
0      & 0      & 0        & 0        & A^2       & -2AN    & 2AL    & 2AC   & L^2-2AR-2CN \\
0      & 0      & 0        & 0        & 0         & 2A^3    & 0      & 0     & 6A^2C \\
0      & 0      & 0        & 0        & 0         & 0       & A^3    & 0     & A^2L-2A^2E \\
0      & 0      & 0        & 0        & 0         & 0       & 0      & A^3   & -A^2D-2A^2N \\
0      & 0      & 0        & 0        & 0         & 0       & 0      & 0     & 2A^4
  \end{pmatrix},
$$
$$
\mathcal M_1=
\begin{pmatrix} 
A      & 0        \\
0      & A      
  \end{pmatrix}
$$
and
$$
\mathfrak D = 4A^{21} \quad\text{and}\quad  \mathfrak D_1=A^2.
$$
\end{proof}

{\sc Remark.\,}
We note that the result remains in force for local algebras over an arbitrary field of the characteristics 0.

\end{document}